\documentclass[review]{elsarticle}

\usepackage{lineno,hyperref}
\modulolinenumbers[5]

\journal{Journal of \LaTeX\ Templates}

\usepackage{algorithm}
\usepackage{algorithmic}
\usepackage{color}
\usepackage{multirow}
\usepackage{supertabular}
\usepackage{longtable}
\usepackage{subfigure}
\usepackage{caption}
\usepackage{amsthm}
\usepackage{amsmath}
\usepackage{amssymb}
\usepackage{hyperref}
\usepackage{url}
\usepackage{makecell}

\begin{document}

\begin{frontmatter}
\title{Generating Linear programming Instances with Controllable Rank and Condition Number}
%\title{Elsevier \LaTeX\ template\tnoteref{mytitlenote}}
%\tnotetext[mytitlenote]{Fully documented templates are available in the elsarticle package on %\href{http://www.ctan.org/tex-archive/macros/latex/contrib/elsarticle}{CTAN}.}


\author[mymainaddress]{Anqi Li}
\ead{lianqi20@mails.ucas.ac.cn}

\author[mymainaddress]{Congying Han\corref{mycorrespondingauthor}}
\cortext[mycorrespondingauthor]{Corresponding author}
\ead{hancy@ucas.ac.cn}

\author[mymainaddress]{Tiande Guo}
\ead{tdguo@ucas.ac.cn}

\address[mymainaddress]{School of Mathematical Sciences, University of Chinese Academy of Sciences, Shijingshan District, Beijing, China}


\begin{abstract}
Instances generation is crucial for linear programming algorithms, which is necessary either to  find the optimal pivot rules by training learning method or to evaluate and verify corresponding algorithms. This study proposes a general framework for designing linear programming instances based on the preset optimal solution. First, we give a constraint matrix generation method with controllable condition number and rank from the perspective of matrix decomposition. Based on the preset optimal solution, the bounded feasible linear programming instance is generated with the right-hand side and objective coefficient satisfying the original and dual feasibility. In addition, we provide three kind of neighborhood exchange operators and prove that instances generated under this method can fill the whole feasible and bounded case space of linear programming. We experimentally validate that  the proposed schedule can generate more controllable linear programming instances, while neighborhood exchange operator can construct more complex instances.
\end{abstract}

\begin{keyword}
%\texttt
{linear programming}\sep {instance generation}\sep {matrix decomposition}\sep {neighborhood exchange operators}
\MSC[2010] 00-01\sep  99-00
\end{keyword}

\end{frontmatter}

\section{Introduction}
The instance generation of linear programming plays an important role in the evaluation and verification of corresponding algorithms. Proper instances can effectively verify the performance of algorithms and find problems in algorithms. Notably, instances with abnormal condition number can effectively test the effectiveness of algorithms in extreme cases. However, the current benchmark~\cite{RW_MPC_2,RW_MPC_15,NetlibArtical} and instance construction algorithms~\cite{RW_MPC_13,RW_MPC_26,RW_MPC_19,compare} cannot generate arbitrarily controllable and full of the entire feasible bounded space instances for important features, such as the condition number and rank of the constraint matrix. Furthermore, it is difficult to design linear programming instances to test the performance of the algorithm.
%Generating instances full of the whole space can effectively evaluate general algorithm performance.

In addition, the lack of efficient training sets is one of the most important reasons that affect the development of the simplex algorithm based on supervised learning~\cite{supervised_pivot_1,supervised_pivot_2,supervised_pivot_3}. There are too few instances in the current benchmark~\cite{NetlibArtical} to construct proper training sets. In recent years, although there are different methods for linear programming instance generation, they cannot meet the needs of specific training sets and test cases. Literature~\cite{compare} has generated instances full of the whole feasible and bounded linear programming space. However, it lacks direct control over the key features such as the condition number and the rank of the constraint matrix. Furthermore, the current methods cannot effectively generate linear programming instances at the level condition number, which cannot meet the needs of general training sets for specific fields. Therefore, we hope to generate corresponding instances for any given condition number and rank, so as to effectively generate representative real data or linear programming problems in some specific fields as training sets, to better design and improve traditional pivot rules based on supervised learning.

\begin{figure*}[t]
\centering
\includegraphics[width=1.00\columnwidth]{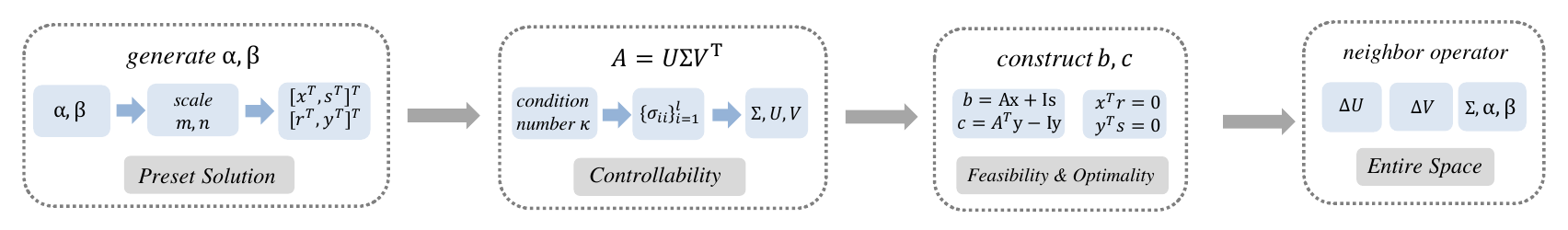}
\caption{Overview of the methodological framework in this paper. The figure consists of four parts: construction of preset solution, the constraint matrix with controllable condition number and rank, the generated right-hand side and objective coefficient satisfying feasibility and optimality, and the neighbor operators designed to fill the entire space.}
\label{fig:overview}
\end{figure*}

Motivated by these observations, we construct a framework to generate instances with controllable condition number and rank, as shown in Figure~\ref{fig:overview}. Specifically, this study focuses on four core aspects: (1) giving an algorithm for constructing preset optimal solution, (2) the constraint matrix $A$ with controllable condition number and rank is presented, (3) the corresponding simplex is constructed to make the preset solution feasible, and the corresponding objective coefficient is designed to make the given solution optimal, (4) three neighborhood exchange operators are designed to ensure that constructed instances can fill the whole feasible and bounded linear programming space.
%First, we propose a method to generate the preset optimal solution. Then the simplex satisfying the feasibility and the objective coefficient guaranteeing the optimality are designed.

Inspired by the literature~\cite{compare}, we construct a set of coding vectors $\alpha,\beta,\alpha \geq 0,\beta_{i} \in \left\{ 0,1 \right\}$ to generate a set of preset optimal solutions of the original and dual linear programming problems. $\alpha$ saves the value corresponding to the basis variables in the original solution and dual solution. $\beta$ is a label vector. $\beta=0$ means that the solution of the corresponding component comes from the dual solution, and $\beta=1$ means that the solution of the corresponding component comes from the original solution. The component in the code corresponds to the component of original solution $X=[x^T, s^T]^T$ and dual solution $Y=[r^T,y^T]^T$. It is worth noting that the derived solution $X=[x^T, s^T]^T, Y=[r^T,y^T]^T$ satisfies $x^{T}r=0, y^T{s}=0$, that is, the complementary relaxation condition. According to the coded dimension, it is natural to derive the size of the generated instance $m=\sum_{i=1}^{N}{\beta_i}$ and $n=N-m$, where $N$ represents the dimension of vector $\beta$, $\beta_i$ is the i-th component of $\beta$, $m,n$ represents the rows and columns of constraints matrix in the original problem.

Based on the preset optimal solution and the size of the generated instance, we need to generate the corresponding simplex as the feasible region. For constraint matrix $A$, we hope to achieve completely control over the condition number and rank. In other words, for any given condition number and rank, we hope to generate the corresponding constraint matrix. Under this goal, the singular value decomposition (SVD) of matrix is a very useful tool. On one hand, SVD can associate the rank with the number of singular values of the constraint matrix. On the other hand, the condition number of the constraint matrix can be effectively expressed as singular values in the sense of two norms. We can perform the inverse process of SVD by generating the corresponding singular value sequence and orthogonal matrices, then generate the constraint matrix with controllable condition number and rank.

Furthermore, we generate the right-hand side to guarantee the feasibility of the preset solution, and the objective coefficient to guarantee its optimality. Considering the standard form of linear programming problem and its dual form, $b$ and $c$ appear on both sides of the equality constraints of the original problem and its dual form. Therefore, we can transform the problem into a standard form by adding relaxation variables, thus deriving the right-hand side and the objective coefficient. Because the solution $X=[x^T, s^T]^T$ and $Y=[r^T,y^T]^T$ are derived from the coding vector groups $\alpha$ and $\beta$, there is only one of $x$ and $r$ takes a non-zero value. $s$ and $y$ are similar. Naturally, the presupposed solution satisfies the complementary relaxation condition. And the constructed simplex guarantees the feasibility and optimality of the preset solution.

Finally, we design neighborhood exchange operators to ensure that the instances generated by neighborhood exchange can fill the entire feasible and bounded linear programming instance space. Therefore, the designed neighborhood exchange operator needs to include all features of the linear programming problem. Considering the SVD of the constraint matrix and $b, c$ derived from the complementary relaxation condition, we design three neighbor operators based on the orthogonal matrix after SVD and the singular value matrix associated with $\alpha$ and $\beta$. In the process of execution, three neighborhood exchange operators can be organized using dynamic adjustment according to probability. In addition, we theoretically prove that instances can fill the entire space under three proposed neighbor operators.

Given exploration on above four aspects, we present a general framework to generate instances with controllable condition number and rank. Based on three neighborhood exchange operators proposed, the generated controllable instances can fill the entire feasible and bounded linear programming instance space. Comprehensive experiments demonstrate that the framework can generate more difficult instances. It is worth noting that compared with the state-of-the-art, our method can generate more complex instances in fewer neighborhood exchanges.

Our main contributions are as follows: 
\begin{itemize}
\item Constructing method of linear programming training set based on supervised learning is proposed.
\item A novel framework for constructing linear programming instances is given based on the preset optimal solution.
\item Provide a constraint matrix generation method with controllable condition number and rank.
\item Three neighbor operators are designed based on SVD, and it is proved that linear programming instances generated based on neighborhood exchange can fill the entire space.
\end{itemize}

%The remainder of this paper is organized as follows. Sections 2 and 3 briefly introduce the background and related works, respectively. In Section 4, we introduce the proposed MCTS rule in detail using four subsections. Section 5 presents experimental results. The conclusions and further work are presented in Sections 6 and 7, respectively.

\section{Related Work}
Test instances is crucial to the algorithm design of linear programming. A good test instance can effectively detect problems existing in the current algorithm, which is convenient for improvement and correction. Constructing controllable instances has received extensive attention, especially for boolean satisfiability~\cite{RW_MPC_1}, quadratic assignment~\cite{RW_MPC_8}, knapsack~\cite{RW_MPC_12, TOU_COR_4}, capacity assignment~\cite{RW_MPC_19}, Hamiltonian completion~\cite{TOU_COR_2} and job shop scheduling~\cite{TOU_COR_5} problem and other instance generation problems~\cite{TOU_COR_1,TOU_COR_3}. These studies show that instance generators are crucial and worth promoting and improving.

Representation and diversity~\cite{compare} are crucial for generating controllable linear programming test instances. However, existing benchmark sets~\cite{RW_MPC_2,RW_MPC_15} do not meet this requirement, as they are usually derived from real-world problem sets or inherited from outdated earlier studies. Although NETLIB~\cite{NetlibArtical} has some linear programming problems that are difficult to solve in reality, the number of instances is too small to meet the design of pivot rules based on supervised learning. MIPLIB~\cite{RW_MPC_20} is a regularly refreshed mixed integer programming test set, but the results are limited by the diversity of submitted questions. Synthetic instances provide another idea for optimization algorithm design. Simple random generation methods tend to construct instances with predictable features~\cite{RW_MPC_1} and cannot meet the requirement of diversity~\cite{RW_MPC_13,RW_MPC_26}. Hooker~\cite{RW_MPC_14} proposes to utilize parametric solvers to construct controllable instances. Although the NETGEN generator~\cite{RW_MPC_19} and its follow-up MNETGEN can generate linear programming instances, they focus on multi-commodity flow, transportation and distribution problems, and cannot generate general controllable feasible linear programming problems.

Todd~\cite{RW_MPC_27} studied the feasible region for generating linear programs using stochastic methods, but it lacked controllability over important features. Pilcher and Rardin~\cite{RW_MPC_24} construct integer programming instances with known partial polytopes based on random cuts. However, such generative methods cannot control instances with interested features, and are limited to specific types of problems like traveling salesman. Smith-Miles~\cite{RW_MPC_26} proposes iterative refinement through search techniques to obtain instances with features of interest. Evolutionary algorithms are commonly used neighborhood search techniques. Hakraborty~\cite{RW_MPC_5} and Cotta~\cite{RW_MPC_6} used this method to make a statistical analysis of the complexity of the algorithm. In addition, there are some works aimed at improving the spectrum of instance hardness~\cite{RW_MPC_28} and diversity of measured features~\cite{RW_MPC_9,RW_MPC_26}. The success of these instance generation methods provides an idea for instance generation for linear programming and integer programming. At the same time, the search technology also provides a starting method for constructing more complete and more complex instance spaces.

\section{Generate Linear Programming Instances}
%In this study, the MCTS rule is proposed to determine the minimum number of pivot iterations and all corresponding pivot paths for the simplex algorithm. This procedure includes the following four parts:
%\begin{itemize}
%\item First, we introduce the proposed SimplexPseudoTree to enable the application of MCTS to the simplex method while avoiding duplicate basic variables.
%\item In the reinforcement learning model stage, we constructed four RL models for the simplex method, which is the basis of applying the Monte Carlo tree search. 
%\item The MCTS rule is proposed and its optimality is proved.
%\item Finally, we find multiple optimal pivot paths for the simplex method and provide a theoretical guarantee.
%\end{itemize}
%The framework is illustrated in Figure~\ref{fig:framework_2}. Next, we describe each part in detail.

\subsection{General Framework for Generating Linear Programming Instances}
\begin{figure*}[t]
\centering
\includegraphics[width=1.00\columnwidth]{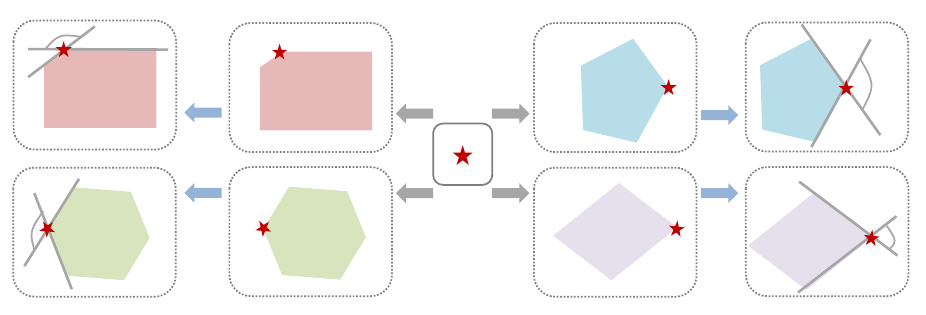}
\caption{A general framework of linear programming instances is constructed. The star represents the preset optimal solution, the gray arrow points to the process of constructing the feasible region, and the blue arrow points to the process of providing the objective coefficient to ensure optimality.}
\label{fig:framework}
\end{figure*}

In the geometric representation of linear programming, simplex represents feasible region and hyperplane represents objective function. The graphical method obtains the point set that maximizes or minimizes the objective function by moving the hyperplane to find the optimal solution. We construct a linear programming problem, specifically, to design the feasible region simplex and the objective function hyperplane. Among all the vertices of the simplex, one must be the optimal solution of the linear programming problem. Therefore, this paper proposes a framework for constructing linear programming instances based on the preset optimal solution, as shown in Figure~\ref{fig:framework}. Specifically, we first generate a point in the space as the preset optimal solution. Then a suitable simplex is constructed as the feasible region. We need to ensure that the preset point is located at the edge of the simplex to ensure that it can be used as the optimal solution of the linear programming instance. Finally, a suitable objective function hyperplane is constructed to ensure the optimality based on the constructed simplex and the preset optimal solution. In the following chapters, we will specifically introduce how to generate each part.

\subsection{Construct Preset Optimal Solution}
\begin{figure*}[t]
\centering
\includegraphics[width=.70\columnwidth]{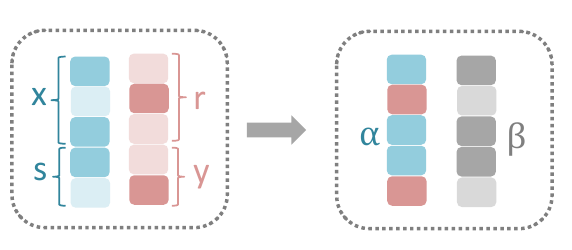}
\caption{Schematic diagram of coding mode.}
\label{fig:AlphaBeta_preset}
\end{figure*}

Inspired by the literature~\cite{compare}, we use the encoding vectors to represent the preset solution. For the original form of linear programming 
\begin{equation}\label{eq:LP_form_primal}
\begin{aligned}
&\max{c^{T}x} \\
&s.t.\ Ax+Is=b \\
&\ \ \ \ \ x\geq0
\end{aligned}
\end{equation}
we record its optimal solution as $X=[x^T,s^T]^T$. Similarly, we record the optimal solution of the dual problem 
\begin{equation}\label{eq:LP_form_dual}
\begin{aligned}
&\min{b^{T}y} \\
&s.t.\ A^{T}y-Iy=c\\
&\ \ \ \ \ y\geq0
\end{aligned}
\end{equation}
as $Y=[r^T,y^T]^T$. The encoding vector consists of $\alpha$ and $\beta$, where the element of each corresponding position of $\alpha$ is the value of the basic variable of the corresponding position in $X$ or $Y$. The optimal solution satisfies $x^{T}r=0$ and $y^{T}s=0$. Therefore, vectors of the same dimension must be able to completely store the values of the basic variables in the original solution and dual solution. $\beta$ is a tag vector where each element is taken as $0$ or $1$. The component with a value of $1$ corresponds to the basic variable of the original problem, and the component with a value of $0$ corresponds to the basic variable of the dual problem, as shown in Figure~\ref{fig:AlphaBeta_preset}.

\begin{algorithm}[H]
\renewcommand{\algorithmicrequire}{\textbf{Input:}}
\renewcommand{\algorithmicensure}{\textbf{Output:}}
\caption{$Encode(a,b,\mu_s,\sigma_s,a^B,\lambda)$}
\label{alg:optimal_coding}
\algsetup{linenosize=\tiny} \scriptsize

\begin{algorithmic}[1]
    \REQUIRE uniform distribution parameter $a$ and $b$, logarithmic normal distribution parameter $\mu_s$ and $\sigma_s$, beta distribution parameter $a^B$, fractional variable proportion parameter $\lambda$.
	\ENSURE coding vectors $\alpha$ and $\beta$, problem scale $m$ and $n$.

    \STATE $N \gets \left\lceil uniform(a,b) \right\rceil$
    \STATE $k \gets \left\lceil uniform(1,N-1) \right\rceil$
    \STATE $I_B \gets $ choose $k$ indices from $\left\{ 1,2,...,N \right\}$ without replacement
    
    \FOR{$i = 1,2,...,N $}
        \STATE $X_1 \gets logNormal(\mu_s,\sigma_s)$
        \STATE $X_2 \gets Beta(a^B,a^B)$
        \IF{$i \in I_B$}
            \STATE $\beta_i=1$
        \ELSE
            \STATE $\beta_i=0$
        \ENDIF 

        \STATE $p_1 \gets uniform(0,1)$    
        \IF{$p_1 \leq \lambda$}
            \STATE $\alpha_i=\left\lceil X_1 \right\rceil-X_2$
        \ELSE
            \STATE $\alpha_i=\left\lceil X_1 \right\rceil$
        \ENDIF
    \ENDFOR
    
    \STATE $m=\sum\nolimits_{i=1}^{N} \beta_i$
    \STATE $n=N-m$
    \end{algorithmic}
\end{algorithm}

We randomly generate the preset optimal coding vector and derive the corresponding coding dimension, as shown in Algorithm~\ref{alg:optimal_coding}. By rounding the random number generated by uniform distribution, we can obtain the total constraint number $N$ of the original problem and dual problem. To ensure that the number of basic variables in the original problem and dual problem is greater than zero, we require that the number of basic variables in the original problem satisfies $1 \leq k \leq N-1$. The basic variable set $I_B$ of the original problem is obtained by sampling without replacement. Scale parameter $\lambda \in (0,1)$ controls whether the components of the optimal basic variables are fractional values. The integer part of the optimal solution is generated by rounding the specified distribution, and the fractional part is generated using the Beta distribution. When the probability random number generated in $(0,1)$ satisfies $p_1 \leq \lambda$, the corresponding component will be taken as a decimal, otherwise it will be taken as an integer. Finally the number of constraints in the original problem is $m=\sum\nolimits_{i=1}^{N} \beta_i$, and the remaining $N-m$ components are the number of constraints in the dual problem.

\subsection{Generating Constraint Matrix with Controllable Condition Number and Rank}
Considering the importance of the condition number to the constraint properties, we hope to generate the corresponding constraint matrix for any given condition number. Notably, the condition number can be transformed into the form of singular values under the $l_2$ matrix norm. Therefore, we need to generate the required singular value sequence for any given condition number, and then construct the constraint matrix with randomness. In addition, singular value decomposition connects the rank of a matrix with its singular values. By controlling the number of singular values, we can further control the rank of the matrix. After generating the singular value sequence and the random orthogonal matrix, we can perform the inverse process of singular value decomposition to obtain the constraint matrix with a given condition number.

\begin{algorithm}[H]
\renewcommand{\algorithmicrequire}{\textbf{Input:}}
\renewcommand{\algorithmicensure}{\textbf{Output:}}
\caption{$generateCond(\kappa,G,a^\kappa,b^\kappa)$}
\label{alg:condition_number}
\algsetup{linenosize=\tiny} \scriptsize

\begin{algorithmic}[1]
    \REQUIRE any given condition number $\kappa$, lower bound $G$ of the difference between the maximum singular value and the minimum singular value, uniform distribution parameter $a^\kappa$ and $b^\kappa$.
	\ENSURE the maximum singular value $\sigma_{min}$, the minimum singular value $\sigma_{max}$.

    \STATE $\sigma_{min} \gets uniform(a^\kappa,b^\kappa)$
    \STATE $\sigma_{max}=\kappa \sigma_{min}$
    
    \WHILE{$\sigma_{max}-\sigma_{min} \leq G$}
        \STATE $\sigma_{min}=2 \sigma_{min}$
        \STATE $\sigma_{max}=\kappa \sigma_{min}$
    \ENDWHILE
    
    \end{algorithmic}
\end{algorithm}

The generation of constraint matrix is shown in Algorithm~\ref{alg:generateA}. Algorithm~\ref{alg:condition_number} is a sub function of Algorithm~\ref{alg:generateA}, which is used to generate the maximum singular value and the minimum singular value for the specified condition number. We first use generateCond function to generate the maximum and minimum non-zero singular values, and then use normal distribution to generate non-zero singular value sequence between maximum and minimum values. The traditional random orthogonal matrix is implemented based on QR decomposition, while the paper~\cite{rvs} improves the problem of relative lack of randomness based on Haar measure. This algorithm has been implemented by the rvs function of python's scipy.stats.ortho\_group package. Finally, we perform the inverse process of singular value decomposition to construct a constraint matrix with controllable condition number and rank.

Theorem~\ref{thm1} proves that the Algorithm~\ref{alg:generateA} can generate the corresponding constraint matrix for any given condition number. Theorem~\ref{thm2} demonstrates that Algorithm~\ref{alg:generateA} can construct the corresponding constraint matrix $A$ for any specified rank.

\newtheorem{thm}{\bf Theorem}[section]
\begin{thm}\label{thm1}
For any $\kappa>0$, Algorithm~\ref{alg:generateA} can construct constraint matrix $A_{\kappa}$ with condition number $\kappa$ under $l_2$ matrix norm.
\end{thm}
\begin{proof}
Considering SVD, the constraint matrix can be reduced to 
\begin{equation}\label{eq:LP_thm1_1}
A_{\kappa}=U_{\kappa} \Sigma_{\kappa} {V_{\kappa}}^T
\end{equation}
where $U_{\kappa} \in \mathbb{R}^{m \times m}$, $V_{\kappa} \in \mathbb{R}^{n \times n}$ and $\Sigma_{\kappa} \in \mathbb{R}^{m \times n}$. $U_{\kappa}$ and $V_{\kappa}$ are orthogonal matrices. The main diagonal elements of $\Sigma_{\kappa}$ are the singular values of constraint matrix $A_{\kappa}$, and the values of other elements are $0$. Notably, the condition number of constraint matrix $A_{\kappa}$ depends on the principal diagonal of $\Sigma_{\kappa}$. Let 
\begin{equation}\label{eq:LP_thm1_2}
\sigma^{\kappa}_{max}=\max{ \left \{ \sigma^{\kappa}_{ii}|i=1,2,...,min(m,n) \right \} }
\end{equation}
be the maximum condition number of the constraint matrix. Similarly, let 
\begin{equation}\label{eq:LP_thm1_3}
\sigma^{\kappa}_{min}=\min{ \left \{ \sigma^{\kappa}_{ii}|i=1,2,...,min(m,n) \right \} }
\end{equation}
be the minimum condition number of the constraint matrix. Under $l_2$ norm, the condition number can be reduced to 
\begin{equation}\label{eq:LP_thm1_4}
\kappa(A_{\kappa})=\left \| A_{\kappa} \right \| \cdot \left \| {A_{\kappa}}^{-1} \right \|=\sigma^{\kappa}_{max} / \sigma^{\kappa}_{min}
\end{equation}
Because the construction of the Algorithm~\ref{alg:generateA} ensures that 
\begin{equation}\label{eq:LP_thm1_4}
\kappa \cdot \sigma^{\kappa}_{min}=\sigma^{\kappa}_{max} \cdot \sigma^{\kappa}_{min}
\end{equation}
Therefore, the condition number under $l_2$ norm is $\kappa$.
\end{proof}

\begin{algorithm}
\renewcommand{\algorithmicrequire}{\textbf{Input:}}
\renewcommand{\algorithmicensure}{\textbf{Output:}}
\caption{$generateA(L,G,a^\kappa,b^\kappa,\kappa,m,n)$}
\label{alg:generateA}
\algsetup{linenosize=\tiny} \scriptsize

\begin{algorithmic}[1]
    \REQUIRE number of nonzero singular values $L$, lower bound $G$ of the difference between the maximum singular value and the minimum singular value, uniform distribution parameter $a^\kappa$ and $b^\kappa$ to generate the maximum singular value and the minimum singular value, preset condition number $\kappa$, number of constraints $m$, number of variables $n$.
	\ENSURE constraint matrix $A$, orthogonal matrix $U$ of $m$ rows and $m$ columns used for SVD, singular value matrix $\Sigma$ used for SVD, orthogonal matrix $V$ of $n$ rows and $n$ columns used for SVD.

    \IF{$L > 1$}
        \STATE $\sigma_{max}, \sigma_{min}=generateCond(\kappa,G,a^\kappa,b^\kappa)$
        \IF{$L > 2$}
            \STATE $\sigma_{vec} \gets $generate $L-2$ values from $uniform(\sigma_{max}, \sigma_{min})$
            \STATE $\sigma_{vec}=sort([\sigma_{vec}, \sigma_{max}, \sigma_{min}])$
        \ELSE
            \STATE $\sigma_{vec}=[\sigma_{max}, \sigma_{min}]$
        \ENDIF    
    \ELSE
        \STATE $\sigma=\kappa$
    \ENDIF
    
    \STATE $\Sigma=zeros(m,n)$
    \STATE $\Sigma(1:L,1:L)=diag(\sigma_{vec})$
    
    \STATE generate random orthogonal matrix $U \in \mathbb{R}^{m \times m}$
    \STATE generate random orthogonal matrix $V \in \mathbb{R}^{n \times n}$
    
    \STATE $A=U \Sigma V^T$
    
    \end{algorithmic}
\end{algorithm}

\begin{thm}\label{thm2}
For any given rank $L \in \mathbb{Z}^{+}$, Algorithm~\ref{alg:generateA} can generate constraint matrix $A_L$ satisfying $rank(A_L)=L$.
\end{thm}
\begin{proof}
Constraint matrix $A_L$ can be written as 
\begin{equation}\label{eq:LP_thm2_1}
A_L=U_L \Sigma_L {V_L}^T
\end{equation}
based on SVD, where $U_L \in \mathbb{R}^{m \times m}$, $V_L \in \mathbb{R}^{n \times n}$ and $\Sigma_L \in \mathbb{R}^{m \times n}$. $U_L$ and $V_L$ are an orthogonal matrices. The main diagonal elements of $\Sigma_L$ are the singular values of constraint matrix $A_L$, and the values of other elements are $0$. Because $U_L$ and $V_L$ are orthogonal matrices, so $U_L$ and $V_L$ must be invertible. Therefore there are a series of elementary matrices $U_1,U_2,...U_u$ and $V_1,V_2,...V_v$ satisfies 
\begin{equation}\label{eq:LP_thm2_1}
U_L=\prod \limits_{i=1}^u{U_i},\ V_L=\prod \limits_{j=1}^v{V_j}
\end{equation}
In this way,
\begin{equation}\label{eq:LP_thm2_3}
A_L=(\prod \limits_{i=1}^u{U_i}) \Sigma_L (\prod \limits_{j=1}^v{V_j})
\end{equation}
Considering that $U_1,U_2,...U_u$ and $V_1,V_2,...V_v$ are elementary matrices, there have
\begin{equation}\label{eq:LP_thm2_3}
rank(A_L)=rank(\Sigma_L)
\end{equation}
Thus the rank of $A_L$ is the number of non-zero singular values of $A_L$. According to Algorithm~\ref{alg:generateA}, $rank(\Sigma_L)=L$, then
\begin{equation}\label{eq:LP_thm2_3}
rank(A_L)=rank(\Sigma_L)=L
\end{equation}
\end{proof}

\begin{algorithm}[H]
\renewcommand{\algorithmicrequire}{\textbf{Input:}}
\renewcommand{\algorithmicensure}{\textbf{Output:}}
\caption{$constructInstance(\alpha,\beta,A,m,n)$}
\label{alg:bc}
\algsetup{linenosize=\tiny} \scriptsize

\begin{algorithmic}[1]
    \REQUIRE coding of preset solution $\alpha,\beta$, constraint matrix $A$, number of constraints $m$, number of variables $n$.
	\ENSURE right-hand side $b$, objective coefficient $c$.

    \FOR{$i = 1,2,...,n+m $}
        \IF{$i \leq n$}
            \STATE $x_i=\beta_i \alpha_i$
            \STATE $r_i=(1-\beta_i) \alpha_i$
        \ELSE
            \STATE $y_{i-n}=(1-\beta_i) \alpha_i$
            \STATE $s_{i-n}=\beta_i \alpha_i$
        \ENDIF
    \ENDFOR
    
    \STATE $b=Ax+Is$
    \STATE $c=A^T y-Ir$

    \end{algorithmic}
\end{algorithm}

\subsection{Generation of the Right-hand Side and the Objective Coefficient}
Considering that the feasible preset solution satisfies $x ^ T r=0$ and $y ^ T s=0$, the preset solution is the optimal solution of the linear programming problem based on the complementary relaxation condition.

The function of the right-hand side is to construct a simplex satisfying the feasibility of the preset solution based on constraint matrix. First, we need to convert the codes corresponding to the preset solution to the original problem $\max c^T x \ \ s.t. Ax+Is=b,x \geq 0$ and the dual problem $\min b^T y \ \ s.t. A^T y+I y=c,y \geq 0$, as shown in the Algorithm~\ref{alg:bc}. The process of constructing original solution and dual solution is the inverse process of coding. The $\alpha_i$ corresponds to $\beta_i=1$ represents the basic variable value of the original problem, and the $\alpha_i$ corresponds to $\beta_i=0$ indicates the value of the basic variable of the dual problem. Thus we can obtain the solution of the problem. Under the given constraint matrix $A$ and the solutions $X=[x ^ T, s ^ T] ^ T$ and $Y=[r ^ T, y ^ T] ^ T$, the right-hand side $b$ and $c$ of the original problem and dual problem can be derived from the linear programming equality constraints, so as to construct a feasible simplex. Considering the feasible preset solution satisfies $x ^ T r=0$ and $y ^ T s=0$, the preset solution is the optimal solution of this linear programming problem based on the complementary relaxation condition.

\begin{algorithm}
\renewcommand{\algorithmicrequire}{\textbf{Input:}}
\renewcommand{\algorithmicensure}{\textbf{Output:}}
\caption{$generate\Sigma(m,n,\Sigma,pos)$}
\label{alg:generateSigma}
\algsetup{linenosize=\tiny} \scriptsize

\begin{algorithmic}[1]
    \REQUIRE number of constraints $m$, number of variables $n$, singular value matrix $\Sigma$, index of neighborhood exchange $pos$.
	\ENSURE updated singular value matrix $\Sigma$.

    \STATE $V \gets$ get nonzero diagnal elements from $\Sigma$
    \STATE $\sigma_{max}=\max(V)$
    \STATE $\sigma_{min}=\min(V)$
    \STATE $\Sigma(pos,pos) \gets uniform(\sigma_{min},\sigma_{max})$
    
    \STATE $V^{'} \gets$ get nonzero diagnal elements from $\Sigma$
    \STATE $V^{'} = sort(V^{'})$
    \STATE $\Sigma \gets$ generate the diagnal maxtrix by $V^{'}$

    \end{algorithmic}
\end{algorithm}

\subsection{Neighborhood Exchange Operators}
\begin{figure*}[t]
\centering
\includegraphics[width=1.00\columnwidth]{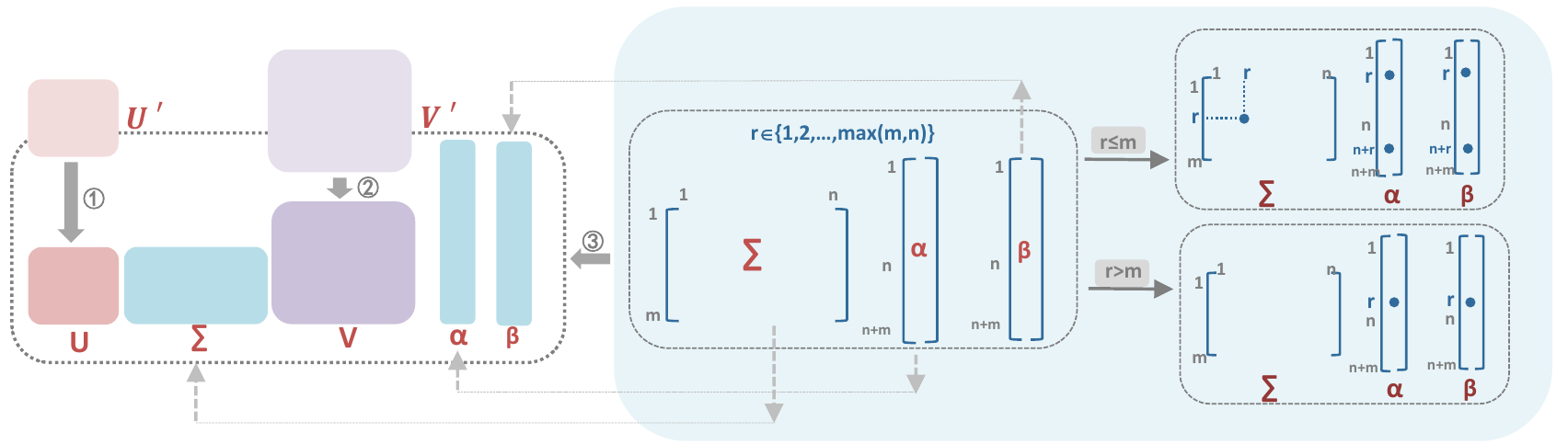}
\caption{Three neighborhood commutative operators. \textcircled{1} represents the first neighbor operator, \textcircled{2} represents the second neighbor operator, and \textcircled{3} represents the third neighbor operator.}
\label{fig:neighbor}
\end{figure*}

The constraint matrix is transformed into two random orthogonal matrices and a matrix with non-zero principal diagonal elements by matrix decomposition. The neighbor operators constructed should also include two random orthogonal matrices of different size, singular values, $\alpha$, $\beta$ and other characteristics. Based on this, we construct three neighborhood exchange operators, as shown in Figure~\ref{fig:neighbor}: (1) generate a random orthogonal matrix $U^{'} \in \mathbb{R}^{m \times m}$ replacing the matrix $U$ in the current SVD; (2) generate a random orthogonal matrix $V^{'} \in \mathbb{R}^{n \times n}$ replacing the matrix $V$ in the current SVD; (3) generate the singular value and $\alpha$, $\beta$ alternative components based on matrix dimension.

\begin{algorithm}[H]
\renewcommand{\algorithmicrequire}{\textbf{Input:}}
\renewcommand{\algorithmicensure}{\textbf{Output:}}
\caption{$generate\alpha\beta(pos,\mu_s,\sigma_s,a^B,\lambda,\xi)$}
\label{alg:generateAlphaBeta}
\algsetup{linenosize=\tiny} \scriptsize

\begin{algorithmic}[1]
    \REQUIRE index of the neighborhood exchange position $pos$, logarithmic normal distribution parameter $\mu_s$ and $\sigma_s$, beta distribution parameter $a^B$, fractional variable proportion parameter $\lambda$, probability $\xi$ of changing a basis variable.
	\ENSURE coding vectors $\alpha$ and $\beta$.

    \STATE $p_1 \gets uniform(0,1)$
    \STATE $X_1 \gets lognormal(\mu_s,\sigma_s)$
    \STATE $X_2 \gets Beta(a^B,a^B)$
    
    \IF{$p_1 \leq \lambda$}
        \STATE $\alpha(pos)=\left\lceil X_1 \right\rceil-X_2$
    \ELSE
        \STATE $\alpha(pos)=\left\lceil X_1 \right\rceil$
    \ENDIF
    
    \IF{$p_1 \leq \xi$}
        \IF {$\beta(pos)=0$}
            \STATE $pos_\beta=\left\{i|\beta_i=1 \right\}$
        \ELSE
            \STATE $pos_\beta=\left\{i|\beta_i=0 \right\}$
        \ENDIF
        
        \STATE $pos_i \gets $rand choose from $pos_\beta$ 
        \STATE $\beta(pos_i)=1-\beta(pos_i)$
        \STATE $\beta(pos)=1-\beta(pos)$
    \ENDIF

    \end{algorithmic}
\end{algorithm}

For the third neighbor operator, we first generate random subscript $r$ in $\left\{1,2,..., max (m, n)\right\}$. Here we suppose that $m \leq n$ (similar for $m > n$). We discuss into two cases: $1 \leq r \leq m$ and $m<r \leq n$. When $1 \leq r \leq m$, the random subscript $r$ corresponds to $r$ in the row subscript and $r$ in the column subscript. In our codes, there are five elements that need to be updated: $\sigma_{rr},\alpha_r,\alpha_{n+r},\beta_r$ and $\beta_{n+r}$. Considering the validity of the updated code, it is necessary to check and correct $\beta$ to satisfy $\sum_{i=1}^{m+n} \beta_i=m$. When $m<r \leq n$, the random subscript $r$ only corresponds to the $r^{th}$ column. In this way, $\alpha_r$ and $\beta_r$ are two elements that need to be updated. It is also necessary to check and correct $\beta$ to satisfy $\sum_{i=1}^{m+n} \beta_i=m$. Algorithm~\ref{alg:condOperator} organizes the above three neighbor operators according to equal probability. In addition, Algorithm~\ref{alg:generateSigma} and Algorithm~\ref{alg:generateAlphaBeta} are two sub functions of Algorithm~\ref{alg:condOperator}, which are used to update the singular value and the preset coding of corresponding positions respectively. Algorithm~\ref{alg:condOperator} can also be organized in other forms depending on our practical needs.

\begin{algorithm}[H]
\renewcommand{\algorithmicrequire}{\textbf{Input:}}
\renewcommand{\algorithmicensure}{\textbf{Output:}}
\caption{$condOperator(m,n,\Sigma,U,V,\lambda,\mu_s,\sigma_s,a_B,\xi)$}
\label{alg:condOperator}
\algsetup{linenosize=\tiny} \scriptsize

\begin{algorithmic}[1]
    \REQUIRE number of constraints $m$, number of variables $n$, the $U,V$ and $\Sigma$ of SVD, fractional variable proportion parameter $\lambda$, logarithmic normal distribution parameter $\mu_s$ and $\sigma_s$, beta distribution parameter $a^B$, probability $\xi$ of changing a basis variable.
	\ENSURE constraint matrix $A$, right-hand side $b$, objective coefficient $c$.

    \STATE $pos \gets $ choose from $\left\{ 1,2,...,\max{(m,n)} \right\}$
    \IF{$pos \leq \max{(m,n)}$}
        \STATE $\Sigma=generate\Sigma(m,n,\Sigma,pos)$
        \STATE $A=U \Sigma V^T$
        \FOR {${pos}^{'} \in \left\{ pos,pos+n \right\}$}
            \STATE $\alpha,\beta=generate\alpha\beta({pos}^{'},\mu_s,\sigma_s,a_B,\lambda,\xi)$
        \ENDFOR
    \ELSIF{$m \leq n$}
        \STATE $\alpha,\beta=generate\alpha\beta(pos,\mu_s,\sigma_s,a_B,\lambda,\xi)$
    \ELSE
        \STATE $\alpha,\beta=generate\alpha\beta(pos+n,\mu_s,\sigma_s,a_B,\lambda,\xi)$
    \ENDIF
    
    \STATE $A,b,c=constructInstance(\alpha,\beta,A,m,n)$

    \end{algorithmic}
\end{algorithm}

\begin{thm}\label{thm3}
Under the above three neighborhood exchange operators, the instances generated by Algorithm~\ref{alg:condOperator} can fill the whole feasible and bounded linear programming instance space.
\end{thm}
\begin{proof}
For any two feasible bounded linear programming instances of the same size, $M$ and $N$. To prove the above conclusion, it is only necessary to show that for instance $M$, instance $N$ can be obtained through finite neighborhood exchanges. In addition, all features of $M$ include $m^M,n^M,A^M,\alpha^M,\beta^M$, and all features of $N$ are $m^N,n^N,A^N,\alpha^N,\beta^N$. Here we suppose that $m \leq n$, and the situation is similar for $m>n$. According to SVD
\begin{equation}\label{eq:LP_thm3_1}
A^M=U^M \Sigma^M (V^M)^T,\ A^N=U^N \Sigma^N (V^N)^T
\end{equation}
First, we use the first neighborhood exchange operator to replace $U^M$ with $U^N$. We perform a neighborhood exchange to obtain 
\begin{equation}\label{eq:LP_thm3_2}
A^M_{(1)}=U^N \Sigma^M (V^M)^T
\end{equation}
Similarly, for the second neighborhood exchange, we use the third neighborhood exchange operator to replace $V^M$ with $V^N$ to obtain 
\begin{equation}\label{eq:LP_thm3_3}
A^M_{(2)}=U^N \Sigma^M (V^N)^T
\end{equation}
Later, we successively take the position and execute the corresponding neighborhood exchange operators for the third kind. Let 
\begin{equation}\label{eq:LP_thm3_4}
\begin{aligned}
&\sigma^M_{rr}=\sigma^N_{rr},\ \ \alpha^M_{r}=\alpha^N_{r},\ \ \alpha^M_{n+r}=\alpha^N_{n+r},\\ &\beta^M_{r}=\beta^N_{r},\ \ \beta^M_{n+r}=\beta^N_{n+r}    
\end{aligned}
\end{equation}
for $r=1,2,...,m$ respectively. In addition, let 
\begin{equation}\label{eq:LP_thm3_5}
\alpha^M_{r}=\alpha^N_{r},\ \ \beta^M_{r}=\beta^N_{r}
\end{equation}
for $r=m+1,m+2,...,n$. In this way, we have 
\begin{equation}\label{eq:LP_thm3_6}
A^M_{(n+2)}=U^N \Sigma^N (V^N)^T,\ \ \alpha^M_{(n+2)}=\alpha^N_r,\ \ \beta^M_{(n+2)}=\beta^N_r
\end{equation}
The instance $M$ is identical to that of $N$.
\end{proof}

\begin{table*}[thpb]
\centering
\caption{Seven distributions and corresponding parameter settings.}
\tiny
\setlength{\tabcolsep}{.5mm}
\resizebox{1.00\columnwidth}{!}{
\begin{tabular}{c|c|c}
    \hline
    \textbf{Distribution} & \textbf{Python Function} & \textbf{Parameter}\\
    \hline
    logarithmic normal & lognormvariate & uniform(-2, 2), uniform(0.1, 10) \\
    \hline
    normal & $\left| normal \right|$ & uniform(-2, 2), uniform(0.1, 10) \\
    \hline
    uniform & uniform & 1, 1000 \\
    \hline
    chi-square & chisquare & $\left[ uniform(1,50) \right]$ \\
    \hline
    F & f & $\left[ uniform(1,50) \right], \left[ uniform(1,50) \right]$ \\
    \hline
    gamma & gamma & uniform(0.05,10), uniform(0.05,10) \\
    \hline
    exponential & exponential & uniform(2,20) \\
    \hline
\end{tabular}}
\label{tab:distribution}
\end{table*}

\section{Experiment}
\subsection{Experiment Setting}
The experiment in this paper compares our instance generation method with the method in the literature~\cite{compare}. We continue to use the parameter setting and experimental mode of the comparison method. To standardize the comparison, we use the scale of literature~\cite{compare} to generate $1000$ linear programming instances with $50$ constraints and $50$ variables to compare the two methods. For the coding that generates the preset optimal solution, we use python's random package to implement and compare values generated by seven distributions. The parameters of seven distributions are shown in Table~\ref{tab:distribution}. Other parameters in the algorithm are shown in Table~\ref{tab:parameter}.

\begin{table*}[thpb]
\centering
\caption{Other parameter settings.}
\tiny
\setlength{\tabcolsep}{.5mm}
\resizebox{.50\columnwidth}{!}{
\begin{tabular}{c|c}
    \hline
    \textbf{Parameter} & \textbf{Value} \\
    \hline
    $N$ & 100 \\
    \hline
    $k$ & 50 \\
    \hline
    $\mu_s$ & uniform(-2, 2) \\
    \hline
    $\sigma_s$ & uniform(0.1, 10) \\
    \hline
    $\lambda$ & 0.7 \\
    \hline
    $a^B$ & lognormvariate(-0.2, 1.8) \\
    \hline
    $\kappa$ & uniform(2, 70000) \\
    \hline
    $\xi$ & 0.5 \\
    \hline
    $a^{\kappa}$ & 0.1 \\
    \hline
    $b^{\kappa}$ & 5 \\
    \hline
    $G$ & 3$\kappa$ \\
    \hline
    $L$ & $\left[ uniform(1,min(m,n)) \right]$ \\
    \hline
\end{tabular}}
\label{tab:parameter}
\end{table*}

\subsection{Distribution Comparison of Preset Coding}
The lognormal distribution of preset codes generated in the literature~\cite{compare} has exponential expectation and variance, which often leads to extreme linear programming instances. Therefore, this paper compares seven commonly used distributions to select the distribution of generated codes. Then, the complete and easy to control optimal solution coding vector can be generated. The seven distributions compared are: logarithmic normal distribution, the absolute value of normal distribution, uniform distribution, chi-square distribution, F-distribution, Gamma distribution and exponential distribution. Table~\ref{tab:b_statistical} and Table~\ref{tab:c_statistical} show the statistical information of the average right-hand side and the average objective coefficient obtained by randomly generating $200$ instances.

\begin{table*}[thpb]
\centering
\caption{Statistical information of the average right-hand side under seven distributions generating optimal coding vectors.}
\tiny
\setlength{\tabcolsep}{.5mm}
\resizebox{1.00\columnwidth}{!}{
\begin{tabular}{c|c|c|c|c|c|c|c}
    \hline
    \textbf{Distribution} & \textbf{Min} & \textbf{Low} & \textbf{Q1} & \textbf{Median} & \textbf{Q3} & \textbf{High} & \textbf{Max} \\
    \hline
    logarithmic normal & -9.54e+13 & -3.33e+05 & -2356.37 & 20.74 & 2.18e+05 &  5.48e+05 & 3.31e+13 \\
    \hline
    normal & -128.94 & -30.34 & -3.93 & 3.28 & 13.68 &  40.10 & 103.58 \\
    \hline
    uniform & -5.36e+04 & -5142.06 & -1181.02 & 201.23 & 1459.68 & 5420.72 & 4.00e+04 \\
    \hline
    chi-square & -2049.59 & -225.56 & -45.54 & 9.92 & 74.48 & 254.50 & 549.77 \\
    \hline
    F & -4.47e+04 & -17.87 & -2.74 & 1.67 & 7.34 & 22.46 & 8352.69 \\
    \hline
    gamma & -1578.49 & -239.71 & -44.56 & 24.39 & 85.54 & 280.69 & 1081.63 \\
    \hline
    exponential & -379.28 & -124.32 & -28.71 & 5.55 & 35.04 & 130.65 & 341.70 \\
    \hline
\end{tabular}}
\label{tab:b_statistical}
\end{table*}

Table~\ref{tab:b_statistical} and Table~\ref{tab:c_statistical} give the statistical information of the average right-hand side and the average objective coefficient constructed by seven distributions, respectively. Statistical information includes the minimum, the low, the first quartile, the median, the third quartile, the high and the maximum. Table~\ref{tab:b_statistical} and Table~\ref{tab:c_statistical} show that the values generated by lognormal distribution are more inclined to generate larger or smaller values. These extreme optimal components have no reasonable significance for practical problems. Normal distribution and exponential distribution tend to generate smaller optimal solution components, which is more limited. Chi-square distribution and gamma distribution can generate slightly larger values, but they are far less extensive than uniform distribution and F-distribution. However, the low, the first quartile, the median, the third quartile and the high of the values generated by the F-distribution are concentrated in the area with very small absolute values. In contrast, the value generated by uniform distribution is the most realistic and easy to control. Therefore, consistent with the expectation and variance, the value constructed by uniform distribution is the properest to control.

\begin{table*}
\centering
\caption{Statistical information of the average objective coefficient under seven distributions to generate optimal coding vectors.}
\tiny
\setlength{\tabcolsep}{.5mm}
\resizebox{1.00\columnwidth}{!}{
\begin{tabular}{c|c|c|c|c|c|c|c}
    \hline
    \textbf{Distribution} & \textbf{Min} & \textbf{Low} & \textbf{Q1} & \textbf{Median} & \textbf{Q3} & \textbf{High} & \textbf{Max} \\
    \hline
    logarithmic normal & -7.64e+12 & -5.28e+05 & -2.11e+05 & -57.75 & 495.34 &  3.17e+05 & 2.01e+12 \\
    \hline
    normal & -152.36 & -36.86 & -11.14 & -1.88 & 6.01 &  31.73 & 115.61 \\
    \hline
    uniform & -3.72e+04 & -4179.11 & -1159.92 & -207.83 & 852.87 & 3872.07 & 2.96e+04 \\
    \hline
    chi-square & -723.45 & -290.59 & -98.55 & -22.24 & 29.47 & 221.50 & 404.15 \\
    \hline
    F & -2962.59 & -20.98 & -6.45 & -1.40 & 3.24 & 17.76 & 714.54 \\
    \hline
    gamma & -702.79 & -291.37 & -78.65 & -10.86 & 63.16 & 275.88 & 1702.79 \\
    \hline
    exponential & -205.78 & -125.42 & -36.39 & -7.07 & 22.96 & 111.98 & 626.89 \\
    \hline
\end{tabular}}
\label{tab:c_statistical}
\end{table*}

\subsection{Feature Comparison}
\begin{figure*}
\centering
\includegraphics[width=.75\columnwidth]{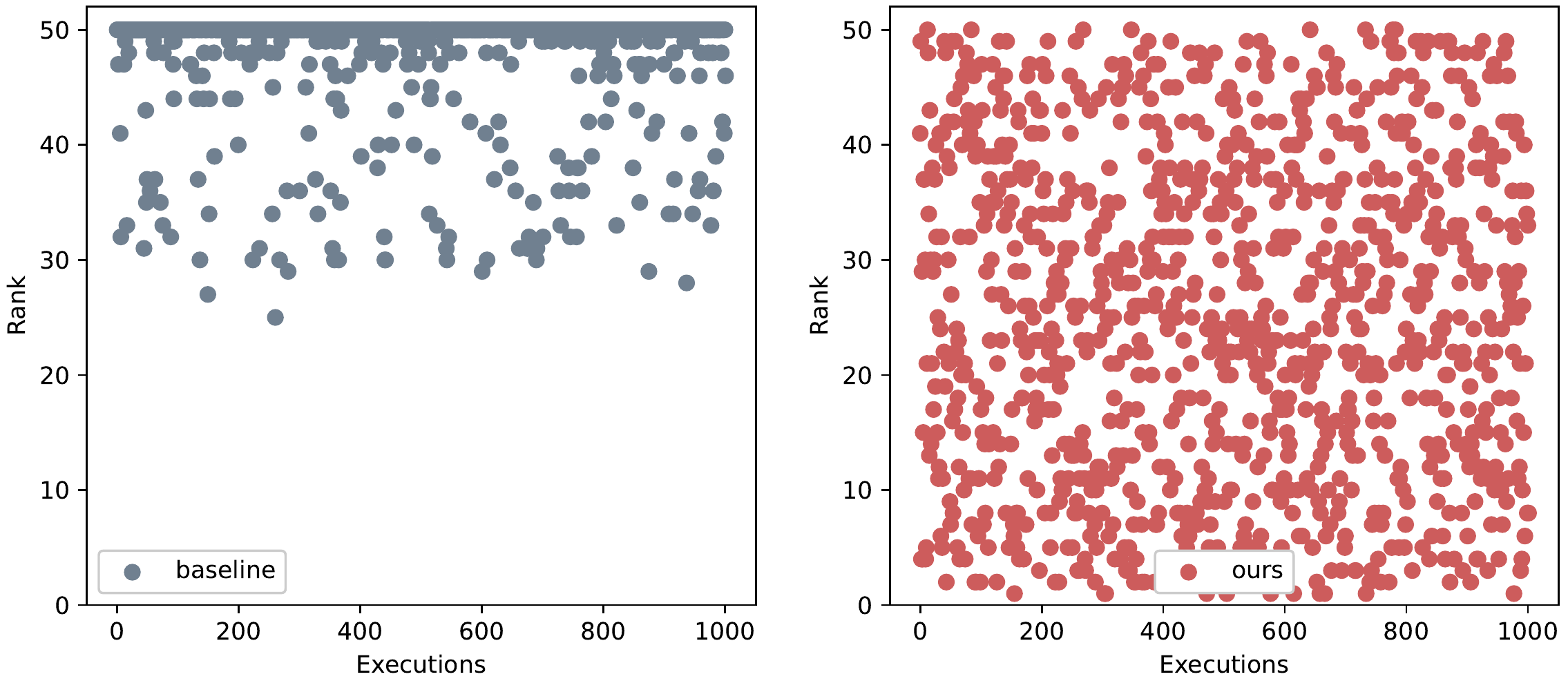}
\caption{Rank distribution of constraint matrix.}
\label{fig:rank}
\end{figure*}

Figure~\ref{fig:rank} shows the distribution diagram of generating $1000$ linear programming instances with $50$ constraints and $50$ variables. Because the rank is less than or equal to the number of constraints, parameters controlling ranks in this algorithm are obtained by rounding the random numbers generated by uniform distribution between $1$ and $50$. Figure~\ref{fig:rank} shows that our method can generate uniform ranks. In contrast, the ranks of the constraint matrix generated by the comparison method are unbalanced. The rank generated within $1000$ executions lacks an interval of $1$ to $25$, which makes it easier to deflect to instances with larger ranks and theoretically uncontrollable.

\begin{figure*}
\centering
\includegraphics[width=.75\columnwidth]{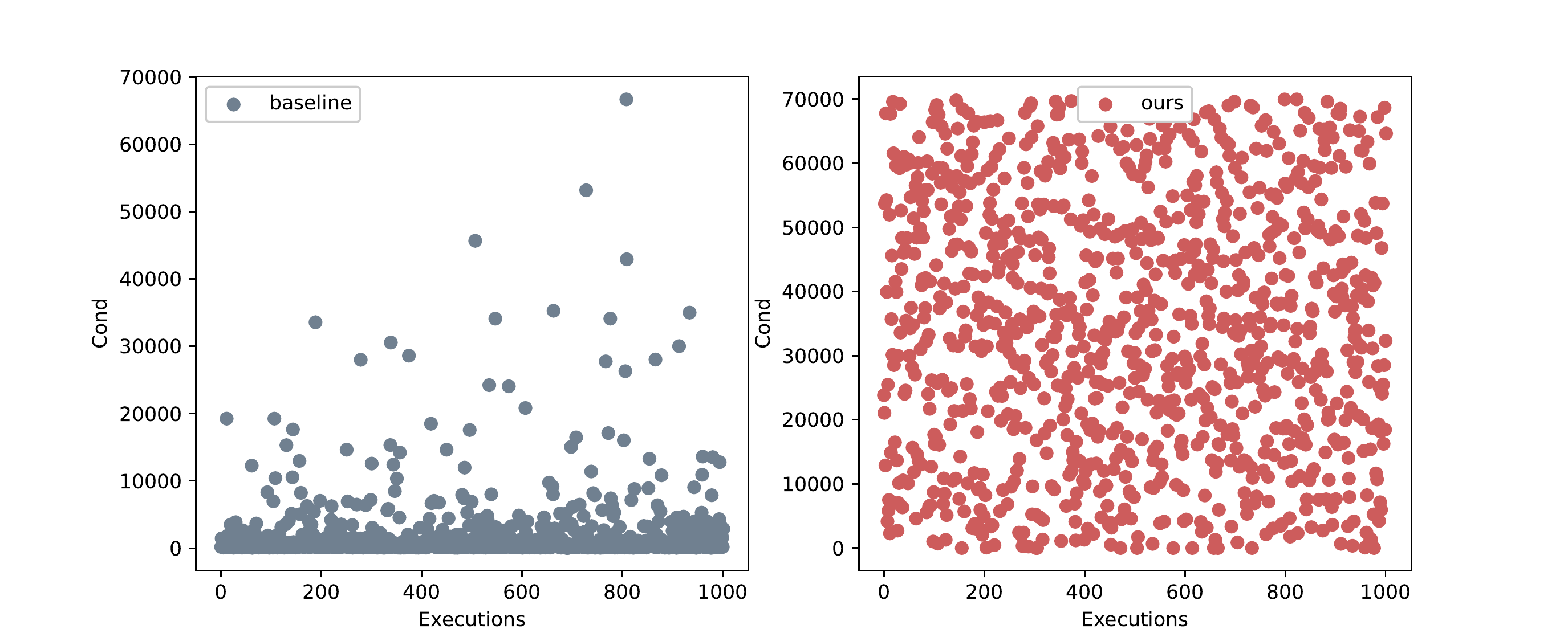}
\caption{The distribution of condition numbers of constraint matrix.}
\label{fig:cond}
\end{figure*}

Figure~\ref{fig:cond} shows the distribution of the condition numbers of the constraint matrices in $1000$ generated linear programming instances with $50$ constraints and $50$ variables. First, we use the literature method~\cite{compare} to generate instances and obtain the distribution of condition numbers. Then, the range of the condition numbers is preset as the corresponding interval. We use uniform distribution to generate the condition numbers, as shown in Figure~\ref{fig:cond}. It can be found that the compared Algorithm~\cite{compare} tend to generate instances with small condition numbers. The generated condition numbers are mainly less than $5000$ and cannot be generated adaptively for the specified values. In contrast, our method can generate a uniformly distributed instance within the specified condition number, which is more suitable for the requirements of linear programming algorithm on test instances.

\subsection{Solver Iterations for Generating Instances}
\begin{figure*}
\centering
\includegraphics[width=1.00\columnwidth]{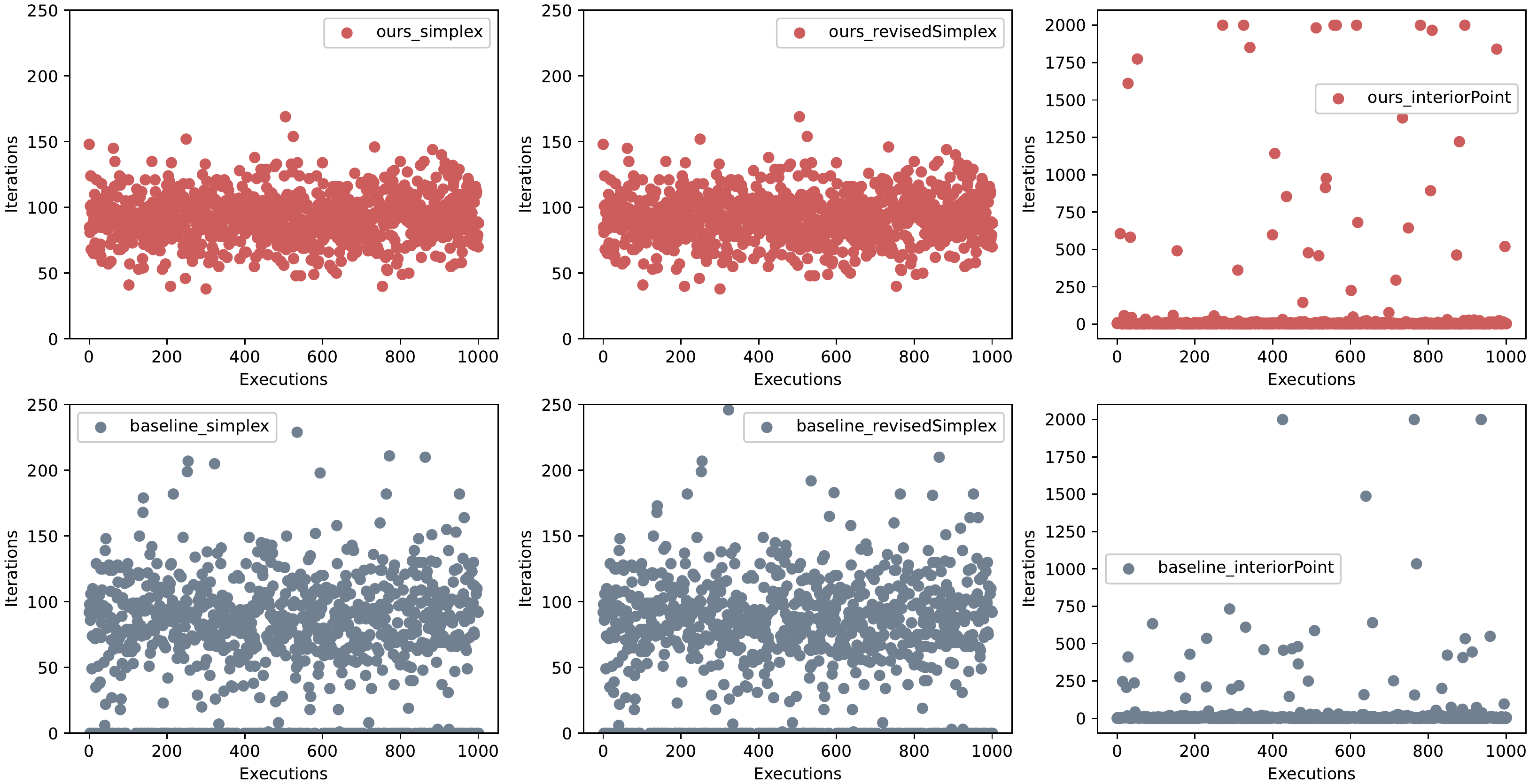}
\caption{Iterations of simplex, revised simplex and interior point method on generated instances.}
\label{fig:linprog_python}
\end{figure*}

Figure~\ref{fig:linprog_python} shows the iterations of three algorithms implemented by the linprog function in python's scipy package on $1000$ generated instances. Ours marks the algorithm proposed in this paper, and baseline marks the comparison method~\cite{compare}. The first column is the number of iterations of $1000$ instances generated by the simplex algorithm in the default configuration. Similarly, the second column corresponds to the revised simplex algorithm, and the third column corresponds to the interior point algorithm. It can be seen from Figure~\ref{fig:linprog_python} that the iterations of our instances executed by the simplex algorithm and the revised simplex algorithm are concentrated in the range of $50$ to $140$. Compared with the literature~\cite{compare}, there are almost no instances with iterations in between $0$ and $50$, which effectively improves the difficulty of generating instances. In addition, the method in the literature~\cite{compare} constructs a few instances with more than $150$ iterations, which can be subsequently implemented by our neighbor operators. For the interior point algorithm, our method is generally similar to the comparison method in the literature~\cite{compare}, but ours is more inclined to generate instances with more than $750$ iterations.

\begin{figure*}[t]
\centering
\includegraphics[width=1.00\columnwidth]{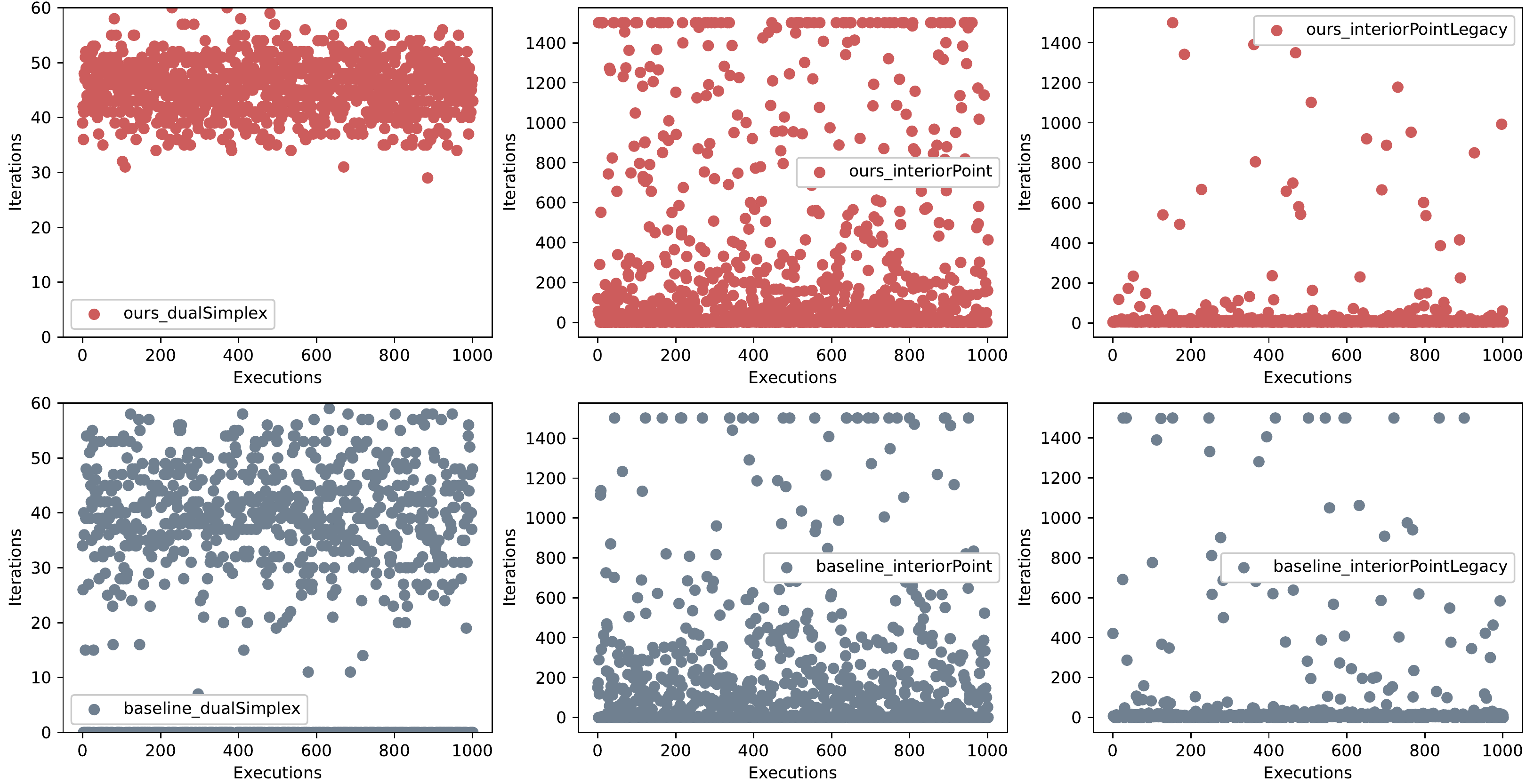}
\caption{Iterations of dual simplex, interior point and interior point legacy method on generated instances.}
\label{fig:linprog_matlab}
\end{figure*}

Figure~\ref{fig:linprog_matlab} shows the iterations of three linear programming algorithms on $1000$ generated instances implemented by matlab's linprog function. Three algorithms are dual simplex algorithm, interior point algorithm and interior point Legacy algorithm, which corresponds to three columns in Figure~\ref{fig:linprog_matlab}. As far as the dual simplex algorithm is concerned, the iterations corresponding to the instances constructed by this method are concentrated in the part above $30$. Compared with the comparison method in the second line, it significantly and effectively improves the difficulty of generated instances. The iteration of the interior point algorithm represented by the second column is similar to that of the python function in Figure~\ref{fig:linprog_python}. In general, there is no difference between our method and the comparison method~\cite{compare}, but ours is more inclined to generate instances with higher iterations. Finally, the two instance generation methods require similar iterations for the inter point Legacy algorithm.

\begin{table*}
\centering
\caption{The comparison results between our neighbor operators and the baseline~\cite{compare}. The first quartile, the median and the third quartile are the statistical information obtained by executing the neighbor operators 100 times. The iteration steps are $0$, $10$, $30$, $50$, $100$, $300$, $500$ and $1000$ respectively.}
\tiny
\resizebox{1.00\columnwidth}{!}{
\begin{tabular}{c|c|ccc|ccc|ccc}
    \hline
    \multicolumn{2}{c}{ } & \multicolumn{3}{c}{Interior Point} & \multicolumn{3}{c}{Simplex} & \multicolumn{3}{c}{Revised Simplex}\\
    \hline
    \textbf{Method} & \textbf{Steps} & \textbf{Q1} & \textbf{Q2} & \textbf{Q3} & \textbf{Q1} & \textbf{Q2} & \textbf{Q3} & \textbf{Q1} & \textbf{Q2} & \textbf{Q3}\\
    \hline
    \hline
    \multirow{8}{*}{\makecell[c]{Ours \\ Search}} & 0 & \textbf{10} & \textbf{10} & \textbf{10} & \textbf{224} & \textbf{224} & \textbf{224} & \textbf{149} & \textbf{149} & \textbf{149} \\
    \cline{3-11}
    &10 & \textbf{10} & \textbf{10} & \textbf{11} & \textbf{226} & \textbf{235} & \textbf{247} & \textbf{192} & \textbf{203} & \textbf{216} \\
    \cline{3-11}
    &30 & \textbf{10} & \textbf{11} & \textbf{11} & \textbf{244} & \textbf{256} & \textbf{269} & \textbf{207} & \textbf{219} & \textbf{227} \\
    \cline{3-11}
    &50 & \textbf{10} & \textbf{11} & \textbf{11} & \textbf{255} & \textbf{266} & \textbf{278} & \textbf{212} & \textbf{224} & \textbf{236} \\
    \cline{3-11}
    &100 & \textbf{11} & \textbf{11} & \textbf{11} & \textbf{264} & \textbf{277} & \textbf{289} & \textbf{222} & \textbf{231} & \textbf{249} \\
    \cline{3-11}
    &300 & \textbf{12} & \textbf{13} & \textbf{14} & \textbf{279} & \textbf{292} & \textbf{314} & \textbf{238} & \textbf{248} & \textbf{261} \\
    \cline{3-11}
    &500 & \textbf{12} & \textbf{13} & \textbf{15} & \textbf{288} & \textbf{301} & \textbf{323} & \textbf{248} & \textbf{260} & \textbf{275} \\
    \cline{3-11}
    &1000 & \textbf{13} & \textbf{14} & \textbf{16} & \textbf{299} & \textbf{319} & \textbf{336} & \textbf{255} & \textbf{269} & \textbf{286} \\
    \hline
    \hline
    \multirow{8}{*}{\makecell[c]{Controllable \\ Search \\ (baseline)}} & 0 & \textbf{10} & \textbf{10} & \textbf{10} & \textbf{224} & \textbf{224} & \textbf{224} & \textbf{149} & \textbf{149} & \textbf{149} \\
    \cline{3-11}
    &10 & \textbf{10} & \textbf{10} & 10 & 224 & 224 & 224 & 173 & 184 & 194 \\
    \cline{3-11}
    &30 & \textbf{10} & 10 & 10 & 244 & 229 & 236 & 188 & 199 & 208 \\
    \cline{3-11}
    &50 & \textbf{10} & 10 & 11 & 226 & 233 & 243 & 194 & 203 & 215 \\
    \cline{3-11}
    &100 & 10 & 10 & 11 & 234 & 245 & 254 & 205 & 217 & 225 \\
    \cline{3-11}
    &300 & 11 & 11 & 12 & 247 & 263 & 275 & 217 & 227 & 240 \\
    \cline{3-11}
    &500 & 11 & 11 & 12 & 254 & 271 & 285 & 223 & 235 & 246 \\
    \cline{3-11}
    &1000 & 11 & 12 & 13 & 267 & 286 & 302 & 232 & 245 & 260 \\
    \hline
\end{tabular}}
\label{tab:neighbor}
\end{table*}

\subsection{Iterative Effect of Neighbor Operators}
\begin{figure*}
\centering
\includegraphics[width=1.00\columnwidth]{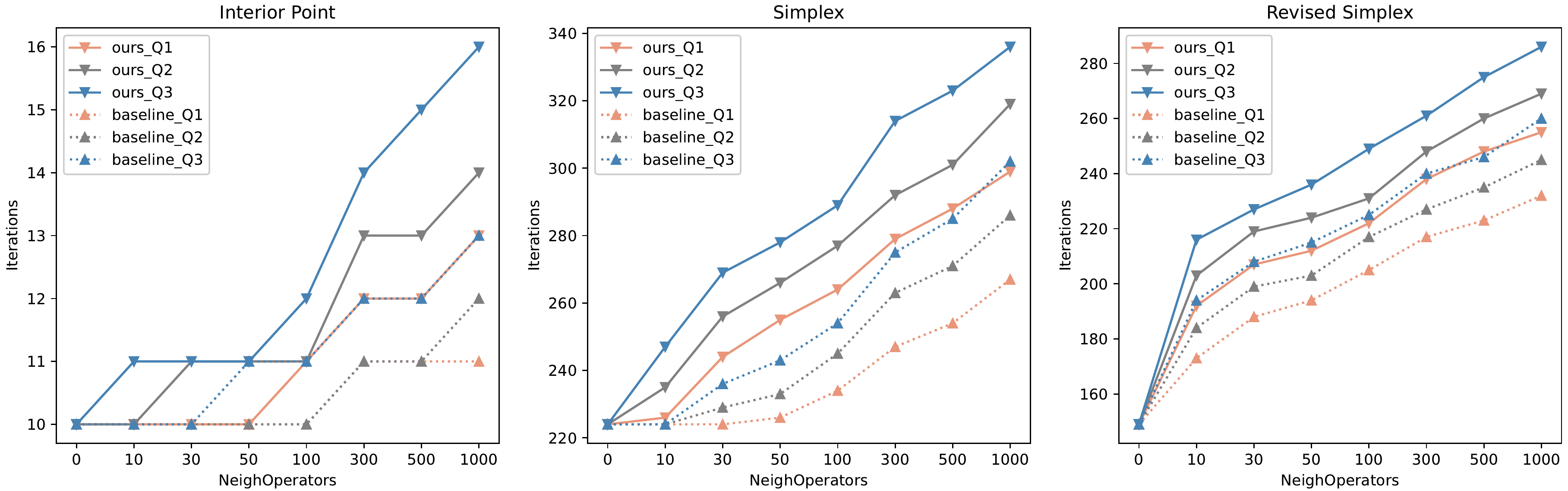} 
\caption{Comparison of our neighbor operators with the baseline~\cite{compare}. The first quartile, the median and the third quartile are the statistical information obtained by executing 100 neighbor operators. The iteration steps are $0$, $10$, $30$, $50$, $100$, $300$, $500$\ and $1000$, respectively. }
\label{fig:neighborOperator_LineChart}
\end{figure*}

We also compare the effect of the neighbor operators proposed in this paper with that in the literature~\cite{compare} on constructing instances with higher iterations. We compare the improvement of the iterations of the inter point algorithm, the simplex algorithm and the reviewed simplex algorithm integrated with the linprog function in python's scipy package by default. Starting from the same instance, we execute each algorithm $100$ times and extract its corresponding statistical information: the first quartile, the median, and the third quartile. Under $0 $, $10 $, $30 $, $50 $, $100 $, $300 $, $500 $and $1000 $ times of neighborhood exchanges, our method can make the inter point algorithm, simplex algorithm, and revised simplex algorithm get higher iterations. To some extent, our neighbor operators are more effective in constructing more complex instances. Figure~\ref{fig:neighborOperator_LineChart} shows the visualization results of Table~\ref{tab:neighbor}. The same color represents the same statistical information, the solid line represents our method, and the dotted line represents the compared method~\cite{compare}. We can find that for curves of the same color, the solid line is always above the dotted line, and the gap is quite large. In other words, for the same iterations, the neighbor operator proposed in this paper can get complexer instances. That means our neighbor operators only need a few iterations to get instances of the same difficulty. In conclusion, ours shows great advantages, compared with the neighbor operators in the literature~\cite{compare}.

\section{Conclusion}
The instance generation framework of linear programming with controllable condition number and rank is proposed from the perspective of preset optimal solution. We first derive the problem size by generating the optimal solution coding. Then the constraint matrix with controllable condition number and rank is constructed using the inverse process of singular value decomposition. In addition, the primal form and dual form of linear programming are used to construct a simplex satisfying the feasibility of the preset solution. At the same time, the right-hand side and the objective coefficient are derived using the complementary relaxation condition. Furthermore, we propose three neighborhood exchange operators to generate more difficult instances. This framework can generate instances with controllable condition number and rank for different linear programming algorithms, and provides a method for constructing training data of simplex algorithm based on supervised learning.

\bibliography{mybibfile}

\end{document}